\documentclass[12pt,a4paper]{article}

\usepackage{amsmath}    
\usepackage{amssymb}
\usepackage{graphicx}   
\usepackage{verbatim}   
\usepackage{color}      

\usepackage{hyperref}   

\usepackage{enumerate}
\usepackage{mathrsfs}
\usepackage{empheq}
\usepackage{bbold}

\usepackage[final]{showlabels}

\usepackage{amsthm}
\newtheorem{definition}{Definition}[section]

\newtheorem{example}{Example}[section]

\newcommand{\del}{\partial}
\renewcommand{\theta}{\vartheta}
\renewcommand{\phi}{\varphi}

\newcommand{\dd}{\mathrm{d}}

\newcommand{\diag}{\mathrm{diag\,}}

\renewcommand{\title}{A hybrid finite element--finite volume method for conservation laws}

\newcommand{\authorOne}{Rémi Abgrall\footnote{Institute for Mathematics and Computational Science, University of Zurich, Winterthurerstrasse 190, 8057 Zurich, Switzerland}}
\newcommand{\authorTwo}{Wasilij Barsukow\footnote{Institute for Mathematics, CNRS UMR 5251, University of Bordeaux, 351 Cours de la Libération, 33405 Talence, France}}

\begin{document}

\begin{center} \Large
\title

\vspace{1cm}

\date{}
\normalsize

\authorOne, \authorTwo
\end{center}

\begin{abstract}

We propose an arbitrarily high-order accurate numerical method for conservation laws that is based on a continuous approximation of the solution. The degrees of freedom are point values at cell interfaces and moments of the solution inside the cell. To lowest ($3^\text{rd}$) order this method reduces to the Active Flux method. The update of the moments is achieved immediately by integrating the conservation law over the cell, integrating by parts and employing the continuity across cell interfaces. We propose two ways how the point values can be updated in time: either by first deriving a semi-discrete method that uses a finite-difference-type formula to approximate the spatial derivative, and integrating this method e.g. with a Runge-Kutta scheme, or by using a characteristics-based update, which is inspired by the original (fully discrete) Active Flux method. We analyze stability and accuracy of the resulting methods.



\end{abstract}

\section{Introduction}

Consider an $m \times m$ system of conservation laws
\begin{align}
 \del_t q + \del_x f(q) &= 0 & q & \colon \mathbb R^+_0 \times \mathbb R \to \mathbb R^m \label{eq:conslaw}\\
 && f &\colon \mathbb R^m \to \mathbb R^m
\end{align}

We shall consider its numerical solution on a grid of cells $[x_{i-\frac12}, x_{i+\frac12}]$, centered at $x_i = i \Delta x$, $i \in I \subset \mathbb Z$. For simplicity, the grid shall be assumed equidistant with $x_{i+\frac12} - x_{i-\frac12} = \Delta x$ $\forall i \in I$. The discrete times at which the numerical solution is obtained are denoted by $(t^n)_{n \in \mathbb Z}$ and the difference $t^{n+1} - t^n$ is generically denoted by $\Delta t$, but the time steps are not assumed to have equal lengths.

Three approaches to the numerical solution of \eqref{eq:conslaw} have inspired this work. Since the seminal work \cite{godunov59difference}, numerical methods are obtained by introducing discontinuities at every cell interface (Riemann problems) and by evolving them over short time intervals (exactly or approximately), see e.g. \cite{toro09} for an introduction to Riemann solvers. The stabilization/upwinding provided by the solutions of Riemann problems has been fruitfully carried over to high order Finite Element methods in \cite{cockburn89} by considering a piecewise continuous approximation space (Discontinuous Galerkin methods). Upon multiplication with a test function and integration by parts, the (a priori multi-valued) flux at a cell interface is replaced by a flux obtained from a Riemann solver. Finally, in \cite{vanleer77,eymann13} the Active Flux method was proposed, which considers a globally continuous and piecewise parabolic reconstruction, and achieves stability/upwinding by solving the initial value problem (IVP) with continuous data that arises at any cell interface. It thus is conceptually close to Godunov's idea, but applies the algorithm of reconstruction--evolution--projection to a continuous choice of reconstruction. In \cite{abgrall22}, several ways of extending this method to arbitrarily high order have been studied. They primarily differ by the choice of the degrees of freedom. In this paper, we aim at showing the relation between one of these choices and both hybrid finite elements and the Discontinuous Galerkin method. We then aim at developing it further by showing different ways how the point values can be updated in time.

The paper is organized as follows. In Section \ref{sec:dof} we present the discretization that we propose, and in Section \ref{sec:interpol} we discuss the corresponding polynomial approximation space. In Section \ref{sec:evolution} we propose two ways of integrating the discrete system forward in time.

We denote by $P^N$ the space of univariate real-valued polynomials of degree at most $N$, and by $C^k$ the class of $k$ times continuously differentiable functions.

\section{Degrees of freedom} \label{sec:dof}

Consider a basis $(b_k)_{k = 1, 2, \ldots}$, $b_k \colon [-\frac{\Delta x}{2}, \frac{\Delta x}{2}] \to \mathbb R$ of some linear space $W \subset C^1$ of differentiable functions. Upon multiplication of \eqref{eq:conslaw} with $b_k(x-x_i)$ and integration by parts one obtains
\begin{align}
 \frac{\dd}{\dd t} \int_{x_{i-\frac12}}^{x_{i+\frac12}} q(t, x) b_k(x - x_i) \, \dd x &+ f(q(t, x_{i+\frac12}))b_k\left(\frac{\Delta x}{2}\right) - f(q(t, x_{i-\frac12}))b_k\left(-\frac{\Delta x}{2}\right) \nonumber \\&- \int_{x_{i-\frac12}}^{x_{i+\frac12}} f(q(t, x)) \del_x b_k(x - x_i) \dd x = 0 
\end{align}

Instead of declaring the expansion coefficients of $q$ with respect to some basis as degrees of freedom, the moments 
\begin{align}
 q^{(k)}_i(t) := A_k \int_{x_{i-\frac12}}^{x_{i+\frac12}} q(t, x) b_k(x - x_i) \,\dd x
\end{align}
themselves shall be the degrees of freedom. This avoids the appearance of any kind of mass matrix. The normalization constant $A_k$ shall be chosen later.

The Discontinuous Galerkin method would replace the flux $f(q(t, x_{i+\frac12}))$ by the flux obtained from the solution of a Riemann problem at $x_{i+\frac12}$. Here, instead, we insist on continuity across cell-interfaces and introduce new independent variables $q_{i+\frac12}(t) \simeq q(t, x_{i+\frac12})$. The evolution of moments therefore simply reads
\begin{align}
 \forall k : \qquad \frac{\dd q^{(k)}_i(t)}{\dd t}  &+ A_k \left( f(q_{i+\frac12}(t)) b_k\left(\frac{\Delta x}{2}\right) -  f(q_{i-\frac12}(t)) b_k\left(-\frac{\Delta x}{2}\right) \right) \nonumber \\&- A_k \int_{x_{i-\frac12}}^{x_{i+\frac12}} f(q(t, x)) \del_x b_k(x - x_i) \dd x = 0  \label{eq:updateofmoments}
\end{align}

The update procedure for the point values $q_{i+\frac12}$ is discussed in Section \ref{sec:evolution}. Upon considering only the $0^\text{th}$ moment (i.e. the average), together with the point values at cell interfaces one obtains the degrees of freedom of Active Flux \cite{vanleer77,eymann13}.

\section{Interpolation} \label{sec:interpol}

We are thus led to a (canonical) hybrid finite element approximation (compare e.g. to Prop. 6.10 (Section 6.3.3, p. 66) in \cite{guermondfinele1}):

\begin{definition} \label{def:finele}
 Consider $(b_\ell)_{\ell \in \mathbb Z}$ a basis of some linear space $W \subset C^1$ of functions and a sequence $(A_\ell)_{\ell \in \mathbb Z}$, $A_\ell \in \mathbb R \backslash \{0\}$. The finite element considered here is the triple $(K, V, \Sigma)$ with
 \begin{enumerate}[i)]
 \item $ K$ the interval $\left[-\frac{\Delta x}{2}, \frac{\Delta x}{2} \right]$,
 \item $ V$ the space $P^N$ of real-valued polynomials on $K$ of degree at most $N \geq 2$, and denote by $V' = hom(V, \mathbb R)$ its dual space,
 \item and $\Sigma \subset  V'$ the set of degrees of freedom spanned by $\{ \sigma_{\frac12},  \sigma_{-\frac12},  \sigma_0,  \sigma_1, \ldots,  \sigma_{N-2}$\} given as
 \begin{align}
   \sigma_{\pm\frac12}(v) &:= v\left(\pm\frac{\Delta x}{2}\right) &  \sigma_\ell(v) &:= A_\ell \int_{-\frac{\Delta x}{2}}^{\frac{\Delta x}{2}} b_\ell(x) v(x) \,\dd x \quad \forall \ell = 0, 1, \ldots, N-2
 \end{align}
 for any $v \in V$.
 \end{enumerate}
\end{definition}

The space $V$ of shape functions can be endowed with a basis $\{ B_\frac12, B_{-\frac12}, B_0, B_1, \ldots B_{N-2} \}$ such that
\begin{align}
 \sigma_r(B_s) = \delta_{rs} \quad \forall r, s \in \left\{ \frac12, -\frac12, 0, 1, \ldots, N-2 \right\}
\end{align}

By abusing notation, we shall denote the extensions of elements of $\Sigma$ to $L^1(K)$ by the same names. Then the interpolation operator $I \colon L^1(K) \to V$ reads
\begin{align}
 I(v)(x) \quad:=\!\!\!\!\!\! \sum_{r \in \left\{ \frac12, -\frac12, 0, 1, \ldots, N-2 \right\}} \!\!\!\!\!\!\sigma_r(v) B_r(x) \quad \forall x \in K
\end{align}

Similarly, we define the reconstruction $R \colon \mathbb R^N \to V$ as
\begin{align}
 R(a_\frac12, a_{-\frac12}, a_0, a_1, \ldots, a_{N-2})(x) \quad:=\!\!\!\!\!\! \sum_{r \in \left\{ \frac12, -\frac12, 0, 1, \ldots, N-2 \right\}} \!\!\!\!\!\!a_r B_r(x) \quad \forall x \in K
\end{align}
The global reconstruction $q_\text{recon} \colon \mathbb R \to \mathbb R^m$ is given by
\begin{align}
 q_\text{recon}\Big|_{x \in \left[x_{i-\frac12}, x_{i+\frac12}\right]} = R(q_{i+\frac12},q_{i-\frac12}, q_i^{(0)}, \ldots, q_i^{(N-2)})
\end{align}
Observe that $q_\text{recon} \in C^0$.

Having extended the domain of the degrees of freedom to $L^1(K)$, one might consider instead of $V$ a more general, non-linear space $\tilde V$, in which case the interpolation operator does not have the form above, as $\tilde V$ would not have a basis. An example of such a space can be found in \cite{barsukow19activeflux} in the context of limiting.

\begin{figure}
 \centering
\includegraphics[width=0.48\textwidth]{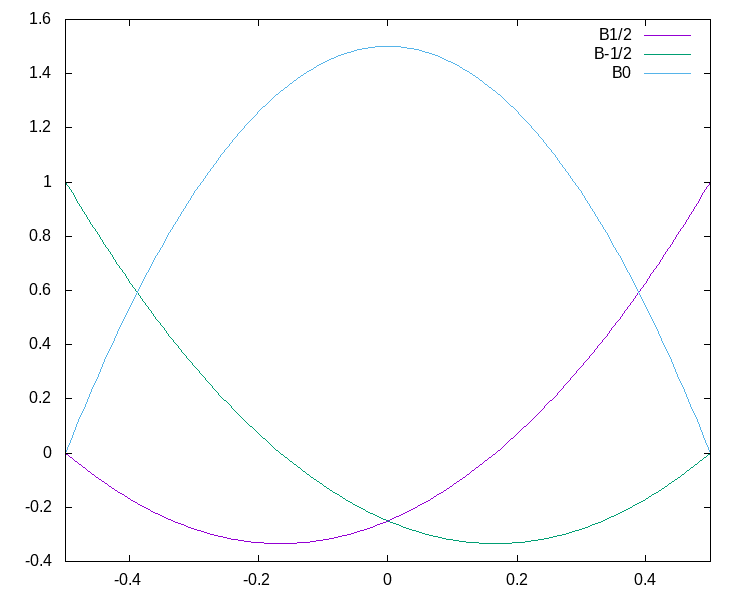} \hfill\includegraphics[width=0.48\textwidth]{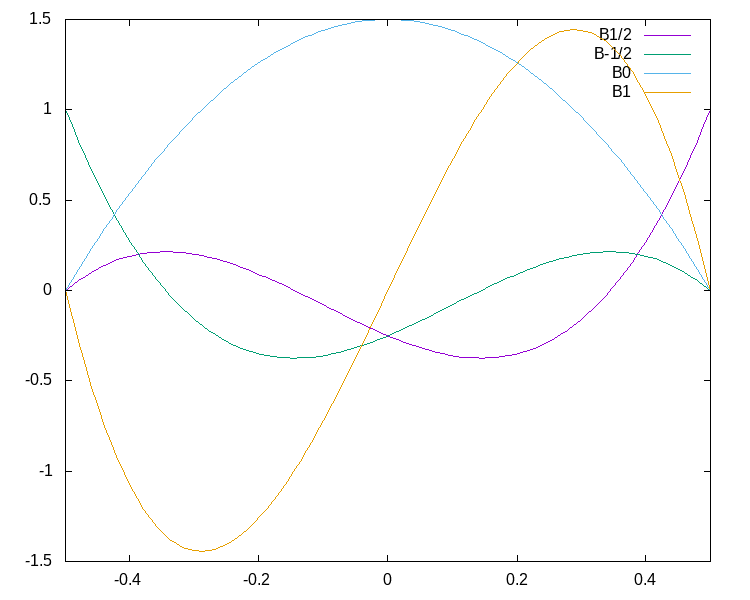} \\
\includegraphics[width=0.48\textwidth]{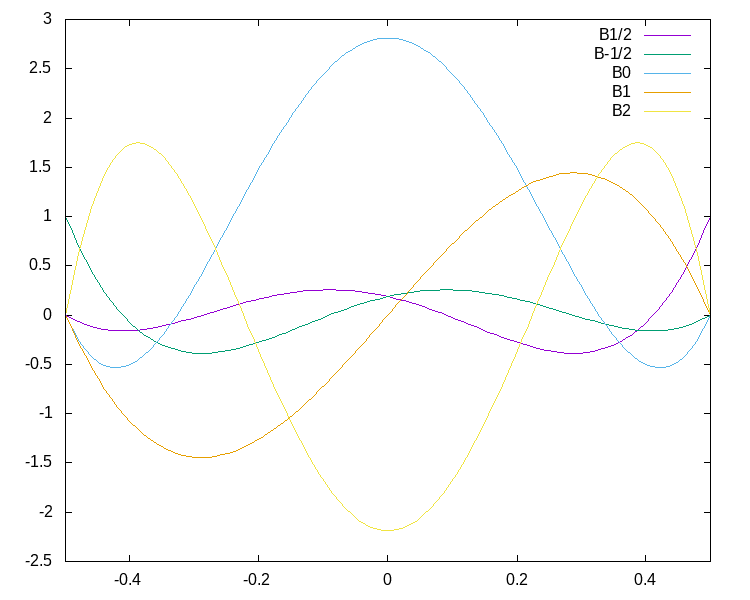} \hfill\includegraphics[width=0.48\textwidth]{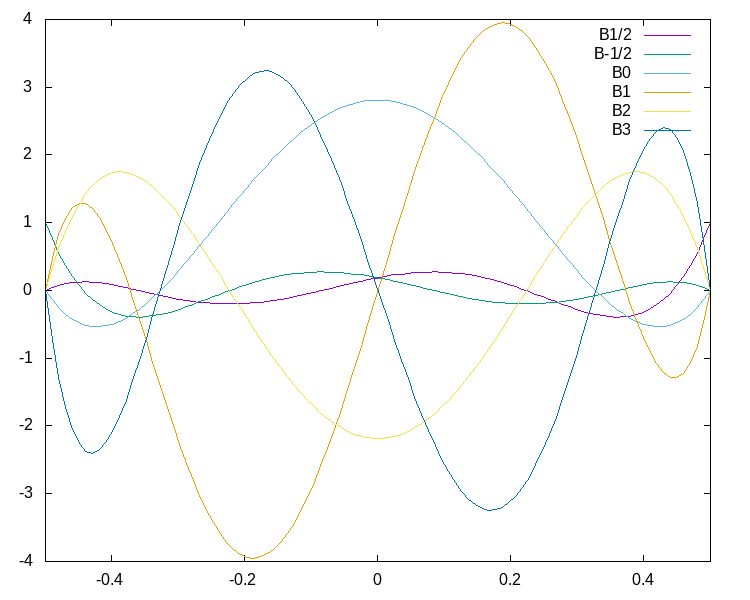} \\
\includegraphics[width=0.48\textwidth]{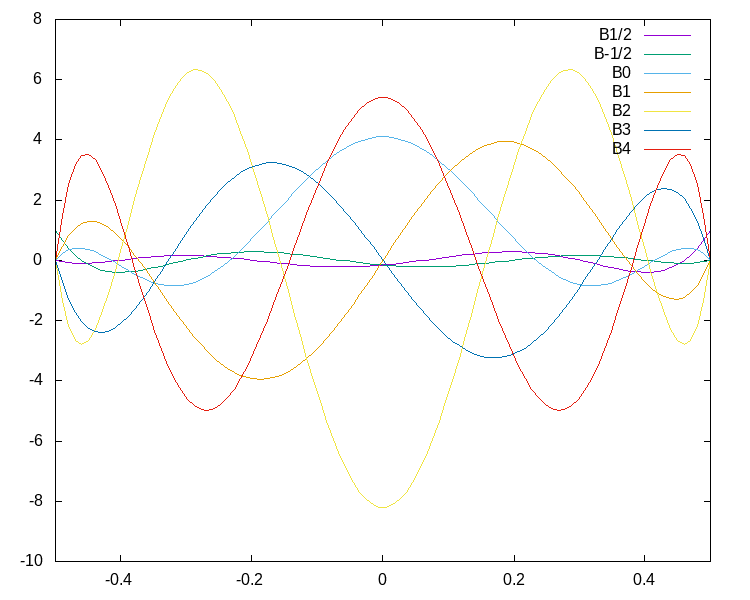}
\caption{The shape functions of Example \ref{ex:monomial}. \emph{Left to right, top to bottom}: $N = 2, \ldots 6$.}
\label{fig:shapefunctions}
\end{figure}

\begin{example} \label{ex:monomial}
 Consider $W = \bigcup_{k \in \mathbb Z} P^k$ with $(b_\ell)_\ell = (1, x, x^2, \ldots)$ its monomial basis. Consider $A_\ell := \frac{(\ell+1) 2^\ell}{\Delta x^{\ell+1}} = \left( \frac{1}{\Delta x}, \frac{4}{\Delta x^2}, \frac{12}{ \Delta x^3} ,\frac{32}{\Delta x^4} , \frac{80}{\Delta x^5} , \ldots \right)$, a choice ensuring that $\sigma_{\ell}(1) = 1$ $\forall \ell \in 2 \mathbb Z$. Then the shape functions (with $\xi := x/\Delta x$) are
 \begin{enumerate}[i)]
  \item for $N=2$
  \begin{align}
B_\frac12 &= \frac{1}{4}  (1+2 \xi) (-1+6 \xi) &
B_{-\frac12} &= \frac{1}{4}  (-1+2 \xi) (1+6 \xi)\\
B_0 &= -\frac{3}{2}  (-1+2 \xi) (1+2 \xi)
  \end{align}
\item for $N=3$
  \begin{align}
B_\frac12 &= \frac{1}{4} (1+2 \xi) \left(-1-4 \xi+20 \xi^2\right) &
B_{-\frac12} &= -\frac{1}{4}  (-1+2 \xi) \left(-1+4 \xi+20 \xi^2\right)\\
B_0 &= -\frac{3}{2}  (-1+2 \xi) (1+2 \xi) &
B_1 &= -\frac{15}{2}  \xi (-1+2 \xi) (1+2 \xi)
  \end{align}
\item for $N=4$
  \begin{align}
B_\frac12 &= \frac{1}{16}  (1+2\xi) \left(3-30 \xi-60 \xi^2+280 \xi^3\right)\\
B_{-\frac12} &= \frac{1}{16}  (-1+2 \xi) \left(-3-30 \xi+60 \xi^2+280 \xi^3\right)\\
B_0 &= \frac{15}{16}  (-1+2 \xi) (1+2 \xi) \left(-3+28 \xi^2\right)\\
B_1 &= -\frac{15}{2}  \xi (-1+2 \xi) (1+2 \xi)\\
B_2 &= -\frac{35}{16} (-1+2 \xi) (1+2 \xi) \left(-1+20 \xi^2\right)
  \end{align}
\item for $N=5$
  \begin{align}
B_\frac12 &= \frac{1}{16} (1+2 \xi) \left(3+24 \xi-168 \xi^2-224 \xi^3+1008 \xi^4\right)\\
B_{-\frac12} &= -\frac{1}{16} (-1+2 \xi) \left(3-24 \xi-168 \xi^2+224 \xi^3+1008 \xi^4\right)\\
B_0 &= \frac{15}{16}  (-1+2 \xi) (1+2 \xi) \left(-3+28 \xi^2\right)\\
B_1 &= \frac{105}{16}  \xi (-1+2 \xi) (1+2 \xi) \left(-5+36 \xi^2\right)\\
B_2 &= -\frac{35}{16}  (-1+2 \xi) (1+2 \xi) \left(-1+20 \xi^2\right)\\
B_3 &= -\frac{315}{32}  \xi (-1+2 \xi) (1+2 \xi) \left(-3+28 \xi^2\right)
  \end{align}
\item for $N=6$
  \begin{align}
B_\frac12 &= \frac{1}{32}  (1+2 \xi) \left(-5+70 \xi+280 \xi^2-1680 \xi^3-1680 \xi^4+7392 \xi^5\right)\\
B_{-\frac12} &= \frac{1}{32}  (-1+2 \xi) \left(5+70 \xi-280 \xi^2-1680 \xi^3+1680 \xi^4+7392 \xi^5\right)\\
B_0 &= -\frac{105}{128}  (-1+2 \xi) (1+2 \xi) \left(5-120 \xi^2+528 \xi^4\right)\\
B_1 &= \frac{105}{16}  \xi (-1+2 \xi) (1+2 \xi) \left(-5+36 \xi^2\right)\\
B_2 &= \frac{105}{64}  (-1+2 \xi) (1+2 \xi) \left(5-232 \xi^2+1232 \xi^4\right)\\
B_3 &= -\frac{315}{32}  \xi (-1+2 \xi) (1+2 \xi) \left(-3+28 \xi^2\right)\\
B_4 &= -\frac{693}{128}  (-1+2 \xi) (1+2 \xi) \left(1-56 \xi^2+336 \xi^4\right)
  \end{align}
They are shown in Figure \ref{fig:shapefunctions}. The evolution equation for the moments becomes
\begin{align}
 \forall k : \qquad \frac{\dd q^{(k)}_i(t)}{\dd t}  &+ (k+1)  \frac{ f(q_{i+\frac12}(t)) -  f(q_{i-\frac12}(t)) (-1)^k }{\Delta x} \nonumber \\&- kA_k \int_{-\Delta x/2}^{\Delta x/2} f(q(t, x + x_i)) x^{k-1} \dd x = 0  \label{eq:updateofmomentsmonomial}
\end{align}
 \end{enumerate}
\end{example}

Observe that for the case of linear flux, the integral in \eqref{eq:updateofmomentsmonomial} can be simplified to
\begin{align}
 kA_k \int_{-\Delta x/2}^{\Delta x/2} f(q(t, x + x_i)) x^{k-1} \dd x &\overset{f(q) = cq}{=} \frac{2(k+1)}{\Delta x} c A_{k-1} \int_{-\Delta x/2}^{\Delta x/2} q(t, x + x_i) x^{k-1} \dd x \\
  &= \frac{2(k+1)}{\Delta x} c q_{i}^{(k-1)} \quad \forall k \geq 1 
\end{align}

For nonlinear flux, we propose to replace this integral by a sufficiently high-order Gauss-Lobatto quadrature of the reconstruction over the cell:
\begin{align}
 &kA_k \int_{-\Delta x/2}^{\Delta x/2} f(q(t, x + x_i)) x^{k-1} \dd x \nonumber \\& \phantom{mmmmmm}\simeq kA_k \Delta x \sum_{m} \omega_m f\Big(R(q_{i+\frac12}(t),q_{i-\frac12}(t),q_i^{(0)}(t),\ldots,q_i^{(N-2)}(t) )(X_m)\Big) X_m^{k-1}  \label{eq:quadratureinspace}
\end{align}
where $(\omega_m)_{m = 0, \ldots, m_\text{max}}$ are the quadrature weights and $X_m \in \left[-\frac{\Delta x}{2},\frac{\Delta x}{2}\right]$ the quadrature points. We use a quadrature formula of $2 N -1$ order of accuracy, such that all moments of polynomial data are integrated exactly.

\section{Evolution}  \label{sec:evolution}

Given that \eqref{eq:updateofmoments} is exact, the properties of the method strongly depend on the choice of the update of the point values $q_{i+\frac12}(t)$. First, in Section \ref{ssec:semidiscrete} a semi-discrete method is presented which is amenable to an integration in time with standard methods, such as SSP-RK. Then, in Section \ref{ssec:fullydiscrete} an alternative method is introduced which is based on an additional integration of \eqref{eq:updateofmoments} over the time step, being thus immediately fully discrete.

\subsection{A method with Runge-Kutta time integration} \label{ssec:semidiscrete}

\begin{definition}[Method A]
Given $(b_k)_k$ and $(A_k)_k$ as in Definition \ref{def:finele} and $N \geq 2$, \emph{Method A} is the semi-discretization
\begin{align}
 \begin{cases} \displaystyle
 \frac{\dd q^{(k)}_i(t)}{\dd t}  =- A_k \left( f(q_{i+\frac12}(t)) b_k\left(\frac{\Delta x}{2}\right) -  f(q_{i-\frac12}(t)) b_k\left(-\frac{\Delta x}{2}\right) \right) \\\displaystyle \phantom{mmmmmmmmmm}+ A_k \int_{-\Delta x/2}^{\Delta x/2} f(q(t, x+x_i)) \del_x b_k(x) \dd x \qquad k = 0, \ldots, N-2 \\ \\
 \displaystyle \frac{\dd}{\dd t} q_{i+\frac12}(t) = -F\Big(q_{i-\frac12}(t),  q_{i}^{(0)}(t), \ldots, q_i^{(N-2)}(t), q_{i+\frac12}(t),  q_{i+1}^{(0)}(t), \ldots, q_{i+1}^{(N-2)}(t), q_{i+\frac32}(t)\Big) 
 \end{cases}
\end{align}
of \eqref{eq:conslaw} with $F$ a consistent approximation of $\del_x f(q)$ at $x_{i+\frac12}$.
\end{definition}

Having in mind the necessity to upwind, i.e. to use information from the left when the characteristic speed $f'$ is positive, and vice versa, we propose to use
\begin{align}
 \frac{\dd}{\dd t} q_{i+\frac12}(t) = - \Big (  f'(\tilde q_{i+\frac12})^+ D  + f'(\tilde q_{i+\frac12})^- D^* \Big ) \label{eq:updatepointvaluemoments}
\end{align}
where, for the moment, one might take $\tilde q_{i+\frac12} = q_{i+\frac12}$. In the scalar case, the positive and negative parts of $f'$ are given by
\begin{align}
 f'(\tilde q_{i+\frac12})^+ &= \max(0, f'(\tilde q_{i+\frac12})) &f'(\tilde q_{i+\frac12})^- &= \min(0, f'(\tilde q_{i+\frac12}))
\end{align}
while for systems, $f'$ is the Jacobian and they are defined through diagonalization:
\begin{align}
 f'(\tilde q_{i+\frac12}) &= T \diag(\lambda_1, \ldots, \lambda_m) T^{-1} & f'(\tilde q_{i+\frac12})^\pm &:= T \diag(\lambda_1^\pm, \ldots, \lambda_m^\pm) T^{-1}
\end{align}
Finally, $D$ and $D^*$ are a left-biased and a right-biased finite difference formula, using the values 
\begin{align}
q_{i-\frac12}(t),  q_{i}^{(0)}(t), \ldots, q_i^{(N-2)}(t), q_{i+\frac12}(t),  q_{i+1}^{(0)}(t), \ldots, q_{i+1}^{(N-2)}(t), q_{i+\frac32}(t)
\end{align}
to approximate $\del_xq(x_{i+\frac12})$. See \cite{iserles82} for high-order finite differences in the context of pure finite difference methods, and compare to \cite{abgrall20}, where those were used after converting moments to point values at cell center.

Arbitrarily high-order finite difference formulae can easily be found by differentiating the reconstruction at the endpoints of $K$:
\begin{align}
 D &:= \frac{\dd}{\dd x} R_i(q_{i+\frac12}, q_{i-\frac12}, q_i^{(0)}, \ldots, q_i^{(N-2)})(\Delta x/2)\\
 D^* &:= \frac{\dd}{\dd x} R_i(q_{i+\frac32}, q_{i+\frac12}, q_{i+1}^{(0)}, \ldots, q_{i+1}^{(N-2)})(-\Delta x/2)
\end{align}

\begin{example} \label{ex:findiff}
 Consider the setup of Example \ref{ex:monomial}. Then, the finite difference formulae read
 \begin{enumerate}[i)]
  \item for $N=2$
  \begin{align}
   \Delta x D &= -6 q_{i}^{(0)}+2 q_{i-\frac12}+4 q_{i+\frac12}\\
   \Delta x D^* &= 6 q_{i+1}^{(0)}-4 q_{i+\frac12}-2 q_{i+\frac32}
   \end{align}
  \item for $N=3$
   
   \begin{align}
\Delta x D &= -6 q_{i}^{(0)}-15 q_{i}^{(1)}-3 q_{i-\frac12}+9 q_{i+\frac12}\\
\Delta x D^* &= 6 q_{i+1}^{(0)}-15 q_{i+1}^{(1)}-9 q_{i+\frac12}+3 q_{i+\frac32}
   \end{align}
  \item for $N=4$
   \begin{align}
\Delta x D &= 15 q_{i}^{(0)}-15 q_{i}^{(1)}-35 q_{i}^{(2)}+4 q_{i-\frac12}+16 q_{i+\frac12}\\
\Delta x D^* &= -15 q_{i+1}^{(0)}-15 q_{i+1}^{(1)}+35 q_{i+1}^{(2)}-16 q_{i+\frac12}-4 q_{i+\frac32}
   \end{align}
  \item for $N=5$
   \begin{align}
\Delta x D &= 15 q_{i}^{(0)}+\frac{105 q_{i}^{(1)}}{2}-35 q_{i}^{(2)}-\frac{315 q_{i}^{(3)}}{4}-5 q_{i-\frac12}+25 q_{i+\frac12}\\
\Delta x D^* &= -15 q_{i+1}^{(0)}+\frac{105 q_{i+1}^{(1)}}{2}+35 q_{i+1}^{(2)}-\frac{315 q_{i+1}^{(3)}}{4}-25 q_{i+\frac12}+5 q_{i+\frac32}
   \end{align}
  \item for $N=6$
   \begin{align}
\Delta x D &= -\frac{105 q_{i}^{(0)}}{4}+\frac{105 q_{i}^{(1)}}{2}+\frac{315 q_{i}^{(2)}}{2}-\frac{315 q_{i}^{(3)}}{4}-\frac{693 q_{i}^{(4)}}{4}+6 q_{i-\frac12}+36q_{i+\frac12}\\
\Delta x D^* &= \frac{105 q_{i+1}^{(0)}}{4}+\frac{105 q_{i+1}^{(1)}}{2}-\frac{315 q_{i+1}^{(2)}}{2}-\frac{315 q_{i+1}^{(3)}}{4}+\frac{693 q_{i+1}^{(4)}}{4}-36 q_{i+\frac12}-6q_{i+\frac32}
  \end{align}
 \end{enumerate}

\end{example}

The method of Example \ref{ex:findiff} has been proposed in \cite{abgrall22} as a high-order extension of the Active Flux method.

After applying the method to scalar linear advection $f(q) = c q$, $q \colon R^+_0 \times \mathbb R \to \mathbb R$, $c > 0$ with the RK3 time integrator we investigate its von Neumann stability by performing the discrete Fourier transform and using the Schur-Miller criterion \cite{schur17,schur18,miller71,barsukow21yee,abgrall22} for the analysis of the eigenvalues of the associated amplification matrix. The maximum CFL numbers $\nu := \frac{\Delta t}{c \Delta x}$ for the setup of Example \ref{ex:findiff} are summarized in Table \ref{tab:cflrk3}. Experimental convergence results are shown in Section \ref{ssec:convergence}.

\begin{table}
 \centering
 \begin{tabular}{c|rrrrr}
  degree $N$ of polynomial  &  2 & 3 & 4 &5 & 6 \\
  order of accuracy         &  3 & 4 & 5 &6 & 7\\
  \hline
  $\text{CFL}_\text{max}$ (Method A, RK3) & 0.41 & 0.21 & 0.13 & 0.09 & 0.06 \\
  $\text{CFL}_\text{max}$ (Method B) & 1.0  & 0.88 & 0.62 & $\simeq 0.6$ & $\simeq 0.6$ 
 \end{tabular}
 \caption{Maximum CFL numbers for linear advection. \emph{Top row}: RK3 time integrator applied to the setup of \ref{ex:findiff}. \emph{Bottom row}: Experimental results for Method B.}
 \label{tab:cflrk3}
\end{table}


\subsection{A method with evolution operators} \label{ssec:fullydiscrete}

As is customary for fully discrete Finite Volume methods, it is also possible to integrate \eqref{eq:updateofmoments} over the time step $[t^n, t^{n+1}]$
\begin{align}
 \forall k : \qquad  \frac{(q^{(k)})^{n+1}_i - (q^{(k)})^{n}_i}{\Delta t}  &+ \frac{A_k}{\Delta t} \int_{t^n}^{t^{n+1}} \left( f(q_{i+\frac12}(t)) b_k\left(\frac{\Delta x}{2}\right) -  f(q_{i-\frac12}(t)) b_k\left(-\frac{\Delta x}{2}\right) \right) \dd t\nonumber \\&- \frac{A_k}{\Delta t} \int_{t^n}^{t^{n+1}} \int_{x_{i-\frac12}}^{x_{i+\frac12}} f(q(t, x)) \del_x b_k(x - x_i) \,\dd x \, \dd t = 0    \label{eq:momentsupdatefullydiscr}
\end{align}
in order to obtain a fully discrete method. It now requires time-averaged fluxes and a space-time quadrature. For $N=2$, the same approach has been followed in the traditional Active Flux method \cite{vanleer77,eymann13}. The point values have been updated using an exact or approximate solution of the initial value problem for \eqref{eq:conslaw} with initial data given by the reconstruction. We thus are led to define

\begin{definition}[Method B]
Given $(b_k)_k$ and $(A_k)_k$ as in Definition \ref{def:finele} and $N \geq 2$, \emph{Method B} is the following discretization of \eqref{eq:conslaw}:
\begin{align}
 \begin{cases} \displaystyle
 (q^{(k)})^{n+1}_i = (q^{(k)})^{n}_i - \Delta t \frac{A_k}{\Delta t} \int_{t^n}^{t^{n+1}} \left( f(q_{i+\frac12}(t)) b_k\left(\frac{\Delta x}{2}\right) -  f(q_{i-\frac12}(t)) b_k\left(-\frac{\Delta x}{2}\right) \right) \dd t \\\displaystyle \phantom{mmmmmmmmmm}+ \Delta t \frac{A_k}{\Delta t} \int_{t^n}^{t^{n+1}} \int_{x_{i-\frac12}}^{x_{i+\frac12}} f(q(t, x)) \del_x b_k(x - x_i) \,\dd x \, \dd t \\ 
%
 \displaystyle q_{i+\frac12}(t) = \Big ( \text{solution at $x=x_{i+\frac12}$ of the IVP \eqref{eq:conslaw} with initial data } q_{\text{recon}} \Big ) + \mathcal O(t^{N+1}) \end{cases}
\end{align}
\end{definition}

For scalar conservation laws and for systems in one spatial dimension this amounts to an estimate of the speed of the characteristic(s):
\begin{itemize}
 \item for linear advection $\del_t q + c \del_x q = 0$:
 \begin{align}
  q_{i+\frac12}(t) = q_\text{recon}(x_{i+\frac12} - c t)
 \end{align}
 \item for Burgers' equation $\del_t q + \del_x (q^2/2) = 0$:
 \begin{align}
  \tilde q_{i+\frac12}^{0} &:= q_\text{recon}(x_{i+\frac12})\\
  \tilde q_{i+\frac12}^{(s)}(t) &:= q_\text{recon}(x_{i+\frac12} - \tilde q_{i+\frac12}^{(s-1)} t) \label{eq:iterationfixpoint}\\
  q_{i+\frac12}(t) &:= \tilde q_{i+\frac12}^{(s_\text{max})}
 \end{align}
 This approximate evolution operator has been derived in \cite{barsukow19activeflux}. Every iteration \eqref{eq:iterationfixpoint} increases the order of accuracy by one (i.e. use $s_\text{max} = N$).
\end{itemize}
For systems in multiple spatial dimensions the solution operator involves characteristic cones. The exact multi-dimensional solution operator for linear acoustics has been used in \cite{eymann13,barsukow18activeflux}.

Method B uses evolution operators, which is close to the original Active Flux method. First, using an exact or approximate evolution operator the following values are computed:
\begin{align}
 q_{i,m}^{n+\ell/\ell_\text{max}} \simeq q(t^{n+\ell/\ell_\text{max}}, x_i + X_m)
\end{align}
for all $\ell = 1, \ldots, \ell_\text{max}$ and as many $m = 0, \ldots, m_\text{max}$ as the Gauss-Lobatto quadrature in space (the same as in \eqref{eq:quadratureinspace}) requires. Note that the times $t^{n+\ell/\ell_\text{max}}$ are not assumed equidistant and the superscript is merely notation; in fact, we propose to choose $\left(t^{n+\ell/\ell_\text{max}}\right)_{\ell = 0,1, \ldots \ell_\text{max}}$ to be Gauss-Lobatto quadrature points in the interval $[t^n, t^{n+1}]$. Observe that all the evaluations of the evolution operator use as initial data the reconstruction at time step $n$. On the one hand, these points allow to perform the space-time quadrature and the time-quadratures in \eqref{eq:momentsupdatefullydiscr}. Some of these newly computed points are kept as new values at cell interfaces:
\begin{align}
 q_{i-\frac12}^{n+1} &:= q_{i,0}^{n+1} & q_{i-\frac12}^{n+1} &:= q_{i,m_\text{max}}^{n+1} 
\end{align}
Observe that the new point values at $t^{n+1}$ are immediately used (through the quadratures) to update the moments to time $t^{n+1}$.

\begin{figure}
 \centering
 \includegraphics[width=\textwidth]{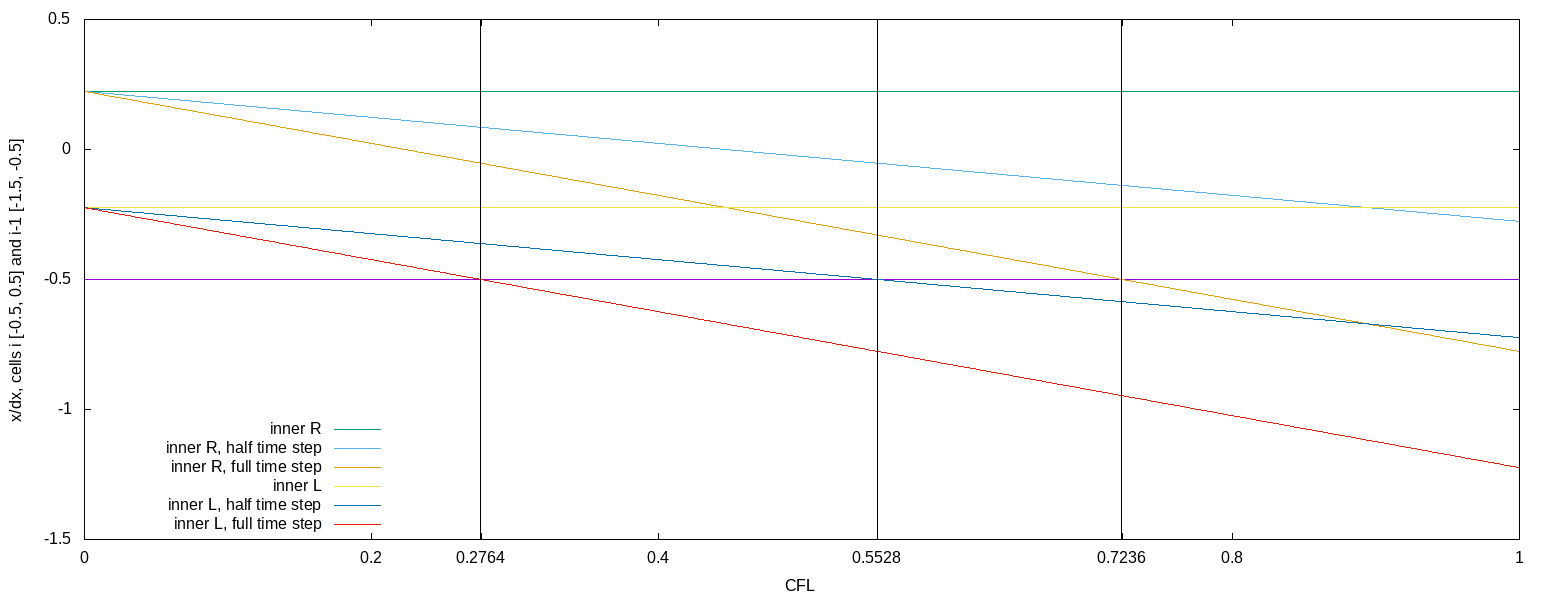} 
 \includegraphics[width=\textwidth]{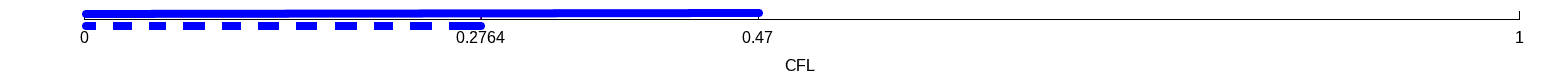} 
 \includegraphics[width=\textwidth]{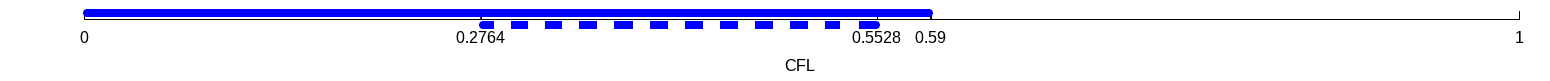} 
 \includegraphics[width=\textwidth]{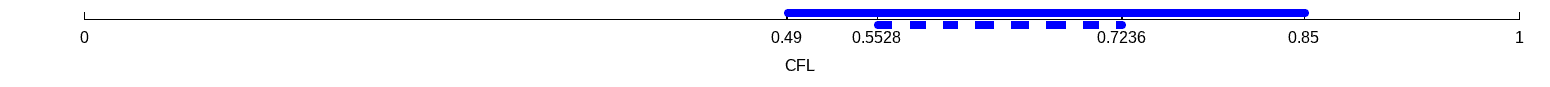} 
 \includegraphics[width=\textwidth]{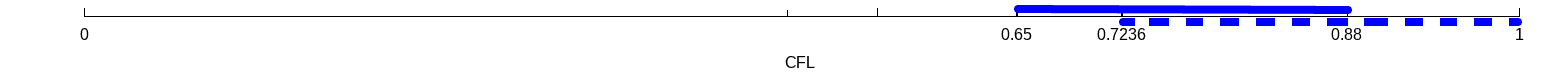} \\
 \includegraphics[width=0.24\textwidth]{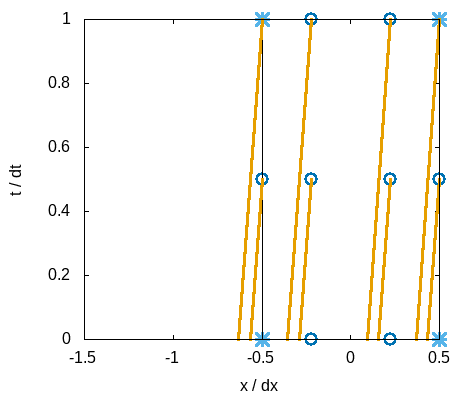} \hfill
 \includegraphics[width=0.24\textwidth]{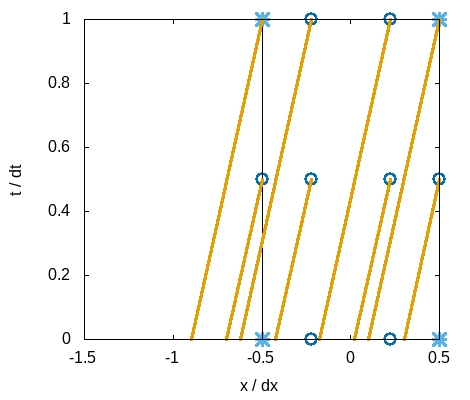} \hfill
 \includegraphics[width=0.24\textwidth]{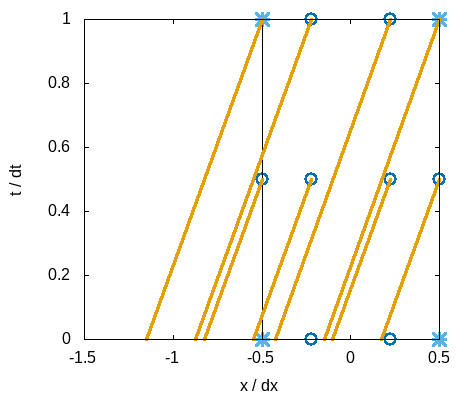} \hfill
 \includegraphics[width=0.24\textwidth]{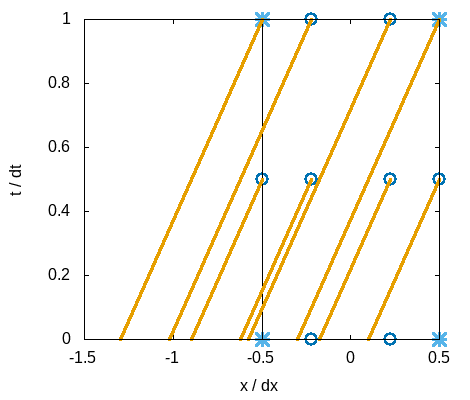} 
 \caption{Method B for $N = 3$: two inner quadrature points in space ($\pm \frac{\Delta x}{2} \frac{1}{\sqrt{5}}$, such that the highest moment ($1^\text{st}$) of the reconstruction is integrated exactly), one inner quadrature point in time (at half time step, Simpson's rule, such that the flux is exact of cubic data)). The Figure shows locations (in the coordinate frame of cell $i$) of the foot points of the characteristics (for linear advection with positive speed) passing through the inner quadrature points at times 0 (trivial), $\frac{\Delta t}{2}$ and $\Delta t$ as a function of the CFL number. Stability analysis requires to determine whether these foot points are located in cell $i$ or $i-1$. This criterion divides the range $[0,1]$ of CFL numbers into four disjoint intervals, as marked by vertical lines. Below, for these four intervals, the physical/expected stability bounds are marked as dashed lines and the actual von Neumann stability as solid lines. The spacetime plots showing the characteristics in the four scenarios are shown as lowest panel. Therein, also the characteristics used for the update of point values at cell interfaces are shown. Points only used for quadrature are shown as circles, while the points retained as independent degrees of freedom are shown as stars.}
 \label{fig:stabilityleapfrog}
\end{figure}

Experimentally, we find this method to be stable with maximal CFL numbers given in Table \ref{tab:cflrk3}, second row, some of which con be confirmed by theoretical analysis. An example of a such for $N=3$ is described in Figure \ref{fig:stabilityleapfrog}. Experimental convergence results are shown in Section \ref{ssec:convergence}.

\begin{figure}
 \centering
 \includegraphics[width=0.49\textwidth]{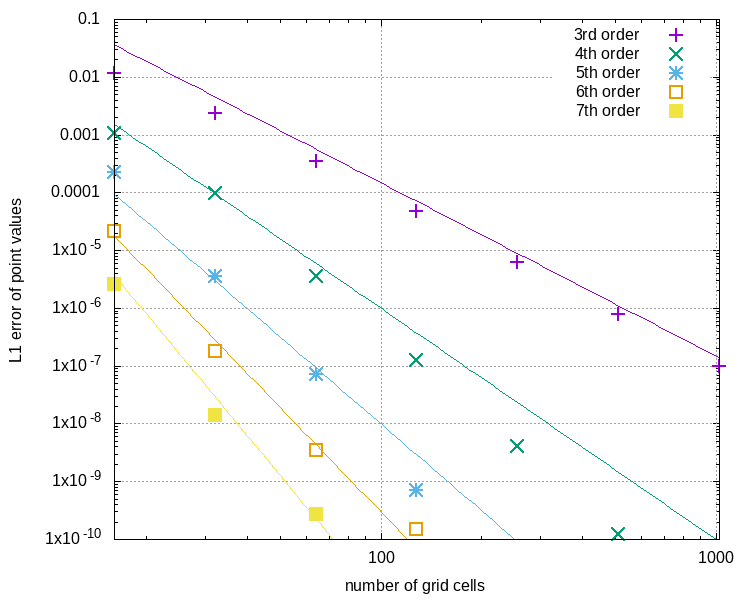} \hfill 
 \includegraphics[width=0.49\textwidth]{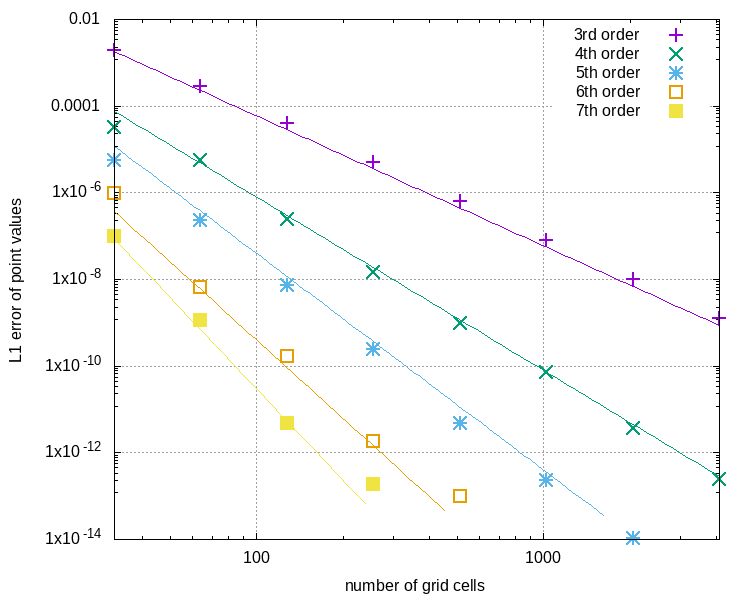} 
 \caption{Convergence analysis of Method A (\emph{left}) and Method B (\emph{right}). Linear advection is solved on grids of various sizes covering the interval $[0,1]$. The initial data are $q_0(x) = 0.8 + \exp(-\frac{(x-\frac12)^2}{0.05^2})$. The $L^1$ error of the point values is shown.}
 \label{fig:errleapfrog}
\end{figure}

\section{Numerical examples}

\subsection{Convergence analysis} \label{ssec:convergence}

We consider linear advection $\del_t q + c \del_x q = 0$ with a speed $c = 1$ on a domain $[0,1]$ endowed with periodic boundaries. For the initial data
\begin{align}
 q_0(x) = 0.8 + \exp\left(-\frac{(x-\frac12)^2}{0.05^2}\right)
\end{align}
a convergence analysis is performed at time $t = 0.1$ for different grids. For Method B, $\text{CFL} = 0.5$ is chosen, and we employ a Runge-Kutta method of third order for Method A. This is why for this latter method, the CFL number is chosen very small ($\text{CFL} = 10^{-4}$), ensuring that the total error is dominated by that of the spatial discretization. The results are shown in Figure \ref{fig:errleapfrog}. The measured orders of convergence agree with the theoretical expectation. One also observes that both methods have a similar magnitude of the error.

\subsection{Burgers' equation}

For nonlinear equations, examples of Method A are shown in \cite{abgrall22}. Method B combines evolution operators and the high-order approach with additional moments. In \cite{barsukow19activeflux}, such approximate evolution operators have been proposed for nonlinear problems; they are endowed with an entropy fix which prevents artefacts (spikes) associated to nearly-transonic shocks. Here, we propose to use the arbitrary-order evolution operator for Burgers' equation from \cite{barsukow19activeflux} and to apply it to the high-degree reconstruction. This evolution operator is finding the correct speed of the characteristic via a fixpoint iteration.

A kind of limiting also has been developed in \cite{barsukow19activeflux} for the third-order Active Flux, which replaces a parabolic reconstruction by a power law when a spurious extremum is recognized. For the details of both the power-law limiting and the entropy fix, the reader is referred to the aforementioned publications. For the high-order method, in \cite{abgrall22} it has been suggested to gradually decrease the order of accuracy if spurious extrema are present. Here, for simplicity, we test at 10 locations inside the cell whether the high-degree polynomial reconstruction exceeds the bounds given by the minimum and maximum of the average and the two point values. In this case, the higher moments are ignored and the parabolic/power-law reconstruction is used, depending on whether the parabola has a spurious extremum.

\begin{figure}
 \centering
 \includegraphics[width=0.49\textwidth]{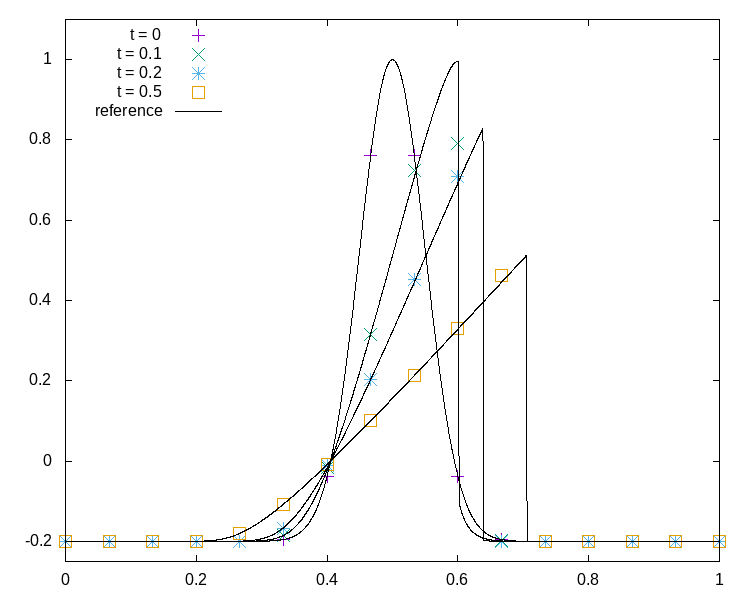} \hfill \includegraphics[width=0.49\textwidth]{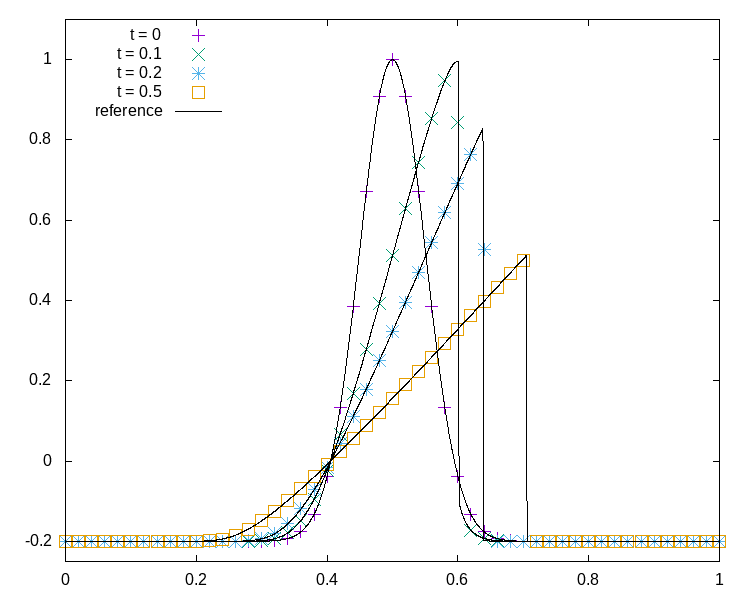}
 \caption{7th order Method B used in conjunction with the iterative evolution operator from \cite{barsukow19activeflux} and limiting as described in the text. Point values are shown. \emph{Left}: Grid of 15 cells. \emph{Right}: Grid of 50 cells.}
 \label{fig:burgers}
 \end{figure}

The self-steepening and shock formation for Gaussian initial data are shown in \ref{fig:burgers} for a 7th order method on very coarse grids. The reference solution has been obtained using the first-order Roe method on a grid of 4000 cells. One observes that the new method has no difficulty to capture the correct solution. Observe also that despite the transonic nature of the setup no artefacts are visible.

\section{Conclusions and outlook}

This paper aims at presenting a new arbitrarily high order numerical method for conservation laws that uses a continuous approximation of the data. It is based on an evolution of the moments of the solution, which is easily obtained through integration of the equation over the cell, and on an independent evolution of point values situated at cell interfaces. This approximation is a high-order extension of the Active Flux method, which is third order accurate, and evolves the point values and only the zeroth moment (cell average). The present method thus can be seen as a hybrid finite element/finite volume method. As the update equations for the moments are exact, many of the details of the method depend on the choice of the evolution of the point values. This paper makes two suggestions: the first is to derive a semi-discrete method by considering finite difference approximations to the flux derivative and to integrate in time using e.g. an SSP-RK method. The second is to use an exact or sufficiently accurate approximate evolution operator based on characteristics. The latter is close to the original approach of the Active Flux method while the former is closer to the ideas presented in \cite{abgrall20}. The latter approach seems to yield a larger stability interval, while the former is significantly easier to apply to complicated problems, such as nonlinear systems of conservation laws. Future work will be devoted to an extension of the proposed method to multiple spatial dimensions.

\end{document}